\UseRawInputEncoding
\documentclass[final,12pt]{elsarticle}
\usepackage{ragged2e}
\justifying
\usepackage{bm}
\usepackage{geometry}
\usepackage{amsmath}
\usepackage{mathrsfs}
\usepackage{amsthm}
\usepackage{amssymb}
\newtheorem{theorem}{Theorem}[section]
\newtheorem{corollary}{Corollary}[section]
\newtheorem{lemma}{Lemma}[section]
\newtheorem{definition}{Definition}[section]

\newenvironment{proof of Theorem 1.1}{\noindent{\textbf{Proof of theorem 1.1.}}\ }{\hfill $\square$\par}

\numberwithin{equation}{section}
\newenvironment{proof1}{\noindent{\textbf{Proof.}}\ }{\hfill $\square$\par}

\usepackage{indentfirst}
\begin{document}
	\begin{frontmatter}  
		\title{Extremal graphs for edge blow-up of lollipops}  
		\author{Yanni Zhai$^1$}  
			\author{Xiying Yuan$^1$\corref{correspondingauthor}}
			\author{Zhenyu Ni$^2$}  
			\cortext[correspondingauthor]{Corresponding author. \\ Email address：{: xiyingyuan@shu.edu.cn} (Xiying Yuan).\\
	        yannizhai2022@163.com(Yanni Zhai).\\
        1051466287@qq.com(Zhenyu Ni).}   
		\address{$^1$ Department of Mathematics, Shanghai University, Shanghai 200444, P.R. China}  
		\address{$^2$ School of Science, Hainan University, Haikou 570228, P.R. China }     
		\begin{abstract}  
			Given a graph $H$ and an integer $p$ ($p\geq 2$), the edge blow-up $H^{p+1}$ of $H$ is the graph obtained from replacing each edge in $H$ by a clique of order $(p+1)$, where the new vertices of the cliques are all distinct. The Tur\'{a}n numbers for edge blow-up of matchings were first studied by Erd\H{o}s and Moon. Very recently some substantial progress of the extremal graphs for $H^{p+1}$ of larger $p$ has been made by Yuan. The range of Tur\'{a}n numbers for edge blow-up of all bipartite graphs when $p\geq 3$ and the exact Tur\'{a}n numbers for edge blow-up of all non-bipartite graphs when $p\geq \chi(H) +1$ has been determined by Yuan (2022), where $\chi(H)$ is the chromatic number of $H$. A lollipop $C_{k,\;\ell}$ is the graph obtained from a cycle $C_k$ by appending a path $P_{\ell+1}$ to one of its vertices. In this paper, we consider the extremal graphs for $C_{k,\;\ell}^{p+1}$ of the rest cases $p=2$ and $p=3$.
			
		\end{abstract}   
	\end{frontmatter}
	\section{Introduction}
	In this paper, we consider undirected graphs without loops and multiedges. For a graph $G$, denote by $E(G)$ the set of edges and $V(G)$ the set of vertices of $G$. The order of a graph is the number of its vertices and the size of a graph is the number of its edges. The number of edges of $G$ is denoted by $e(G)=\lvert E(G)\rvert$. For a vertex $v$ of graph $G$, the neighborhood of $v$ in $G$ is denoted by $N_{G}(v)=\{u\in V(G):  uv\in E(G)\}$. The degree of the vertex $v$, written as $d_{G}(v)$ or simply $d(v)$, is the number of edges incident with $v$. Usually, a path of order $n$ is denoted by $P_{n}$, a star of order $n+1$ is denoted by $S_n$ and a cycle of order $n$ is denoted by $C_{n}$. An independent set of order $n$ is denoted by $I_n$. A matching in $G$ is a set of vertex disjoint edges from $E(G)$, denoted  by $M_k$ a matching of size $k$. For $U\subseteq V(G)$, let $G[U]$ be the subgraph of $G$ induced by $U$, $G-U$ be the graph obtained by deleting all the vertices in $U$ and their incident edges. The graph $K_p(i_1,i_2,\cdots,i_p)$ denote the complete $p$-partite graph with parts of order $i_1,\,i_2,\,\cdots,\,i_p$. Denoted by $T_p(n)$, the $p$-partite Tur\'{a}n graph is the complete $p$-partite graph on $n$ vertices with the order of each partite set as equal as possible. 
	
	For two graphs $G$ and $H$, the union of graphs $G$ and $H$ is the graph $G\cup H$ with vertex set $V(G)\cup V(H)$ and edge set $E(G)\cup E(H)$. In particular, $G=kH$ is the vertex-disjoint union of $k$ copies of $H$. The join of $G$ and $H$, denoted by $G\vee H$, is the graph obtained from $G\cup H$  by adding all edges between $V(G)$ and $V(H)$. Let $H(n,\,p,\,q)$ be the graph $K_{q-1}\vee T_p(n-q+1)$,  and $H'(n,\,p,\,q)$ be any of the graphs obtained by putting one extra edge in any class of $T_p(n-q+1)$. Let $H^{*}(n)$ be graphs obtained by putting (almost) perfect matchings in both classes in $K_2({\lceil \frac{n}{2}\rceil ,\;\lfloor \frac{n}{2} \rfloor }$).
	
	A graph is $H$-free if it does not contain a copy of $H$ as a subgraph. The Tur\'{a}n number ex$(n,\,H)$ is the maximum number of edges in a graph of order $n$ which is $H$-free. Denote by EX$(n,\,H)$ the set of $H$-free graphs of order $n$ with $ex(n,\,H)$ edges and call a graph in EX$(n,\,H)$ an extremal graph for $H$.
	
	Given a graph $H$, the blow-up of $H$, denoted as $H^{p+1}$, is obtained from $H$ by replacing each edge in $H$ by a clique of order $(p+1)$, where the new vertices of the cliques are all different. In 1959, Erd\H{o}s and Gallai \cite{ref2} characterized the extremal graphs for $M_{2k}$. Later, Erd\H{o}s \cite{ref4} studied the Tur\'{a}n numbers of $M_{2k}^3$. Moon \cite{ref9} and Simonovits \cite{ref11} determined the extremal graphs for $M_{2k}^{p+1}$ when $p\geq 2$. Erd\H{o}s, F$\ddot{u}$redi, Gould and Gunderson \cite{ref5} determined the Tur\'{a}n number of $S_{k+1}^3$. Chen, Gould, Pfender and Wei \cite{ref1} determined the Tur\'{a}n number of $S_{k+1}^{p+1}$ for general $p\geq 3$. Glebov \cite{ref6} determined the extremal graphs for edge blow-up of paths. Later, Liu \cite{ref8} generalized Glebov's result to edge blow-up of paths, cycles and a class of trees. Very recently, Wang, Hou, Liu and Ma \cite{ref14} determined the Tur\'{a}n numbers for edge blow-up of a large family of trees. 	A keyring is a $(k+s)$-edge graph obtained from a cycle of order $k$ by appending $s$ leaves to one of its vertices. In some way it is a generalization of cycle and star. 
	Ni, Kang and Shan \cite{ref10} determined the extremal graphs for edge blow-up of keyrings.
A lollipop $C_{k,\;\ell}$ is the graph obtained from a cycle of order $k$ by appending a path $P_{\ell +1}$ to one of its vertices, and the vertex $v\in V(C_{k,\;\ell})$ of degree 3 is called the center of the lollipop. In this paper we will consider the extremal graphs for edge blow-up of lollipops.
	The range of Tur\'{a}n numbers for edge blow-up of all bipartite graphs when $p\geq 3$ and the exact Tur\'{a}n numbers for edge blow-up of all non-bipartite graphs when $p\geq \chi(H) +1$ has been determined by Yuan in \cite{ref15}. The following Theorem 1.1 (i) and (ii) are implied by Theorem 2.3 in \cite{ref15} by taking $\mathcal{B}=\{K_t\}$ when $\ell$ is odd, and $\mathcal{B}=\{K_{t+1}\}$ when $\ell$ is even; and Theorem 1.1 (iii) and (iv) are implied by Theorem 2.4 in \cite{ref15} by taking $\mathcal{B}=\{K_t\}$ when $\ell$ is odd, and $\mathcal{B}=\{K_{t+1}\}$ when $\ell$ is even.

\begin{theorem} \cite{ref15}
	 Suppose $k\geq 3$, $\ell \geq 2$, $t=\lfloor \frac{k-1}{2}\rfloor +\lfloor \frac{\ell -1}{2}\rfloor$, $n$ is sufficiently large.
	 \begin{enumerate}[(i)]
	 	\item When $k$ is even, $\ell$ is odd, $p\geq 3$, ${\rm EX}(n,\,C_{k,\;\ell}^{p+1})=H^{'}(n,\,p,\,t+1)$.
	 	
	 	\item When $k$ is even, $\ell$ is even, $p\geq 3$, ${\rm EX}(n,\,C_{k,\;\ell}^{p+1})=H(n,\,p,\,t+2)$.
	 	
	 	\item When $k$ is odd, $\ell$ is odd, $p\geq 4$, ${\rm EX}(n,\,C_{k,\;\ell}^{p+1})=H(n,\,p,\,t+1)$.
	 
	 \item When $k$ is odd, $\ell$ is even, $p\geq 4$, ${\rm EX}(n,\,C_{k,\;\ell}^{p+1})=H(n,\,p,\,t+2)$.
	 \end{enumerate}
\end{theorem}
	
		How about the extremal graphs for $C_{k,\;\ell}^{p+1}$ for small $p$?  In this paper we will consider the rest cases of Theorem 1.1. More precisely, the extremal graph for $C_{k,\;\ell}^4$ when $k\geq 3$ is odd, $\ell \geq 2$ is determined in Theorem 3.1; the extremal graph for $C_{k,\;\ell}^3$ when $k\geq 4$, $\ell \geq 2$ is odd is determined in Theorem 3.2; the extremal graph for $C_{k,\;\ell}^3$ when $k\geq 4$, $\ell \geq 2$ is even is determined in Theorem 3.3; the extremal graph for $C_{3,\;\ell}^3$ when $\ell \geq 2$  is determined in Theorem 3.4.

		Simonovits proposed the following problem.
		
		\noindent\textbf{Problem 1.1.}\; \cite{ref8} Characterize graphs whose unique extremal graph is of the form $H(n,\,p,\,q)$, where $q\geq 1$, $p\geq 2$.
		
		Combining Theorems 1.1, 1.2, 3.1 - 3.4, the extremal graph for $C_{k,\;\ell}^{p+1}$ is of the form $H(n,\,p,\,q)$. In this way an additional family of forbidden graphs of Problem 1.1 is provided in this paper.
		
		We would like to point out that in \cite{ref10} Ni, Kang and Shan have determined the extremal graphs for $C_{k,\;1}^{p+1}$ (see Theorem 1.2). In this paper we suppose $\ell \geq 2$.
			\begin{theorem}(\cite{ref10})
			When $k\geq 3$, $p\geq 2$, $n$ is sufficiently large, let $G$ be the extremal garph for $C_{k,\,1}^{p+1}$.
			\begin{enumerate}[(i)]
				\item When $(p,k)\neq (2,3)$, $H(n,\,p,\,\lfloor \frac{k-1}{2}\rfloor +1)$ ($H'(n,\,p,\,\lfloor \frac{k-1}{2}\rfloor +1)$ resp.) is the unique extremal graph for $C_{k,\,1}^{p+1}$ when $k$ is odd (even resp.).
				\item When $(p,k)= (2,3)$
				\begin{align*}
					\begin{split}	
						G \in \left\{
						\begin{array}{ll}
							\{(\frac{1}{3} \lceil \frac{n}{2} \rceil K_{3})\vee I_{\lfloor \frac{n}{2} \rfloor },\, H^{*}(n)\}, & {\rm if}\quad 12\vert n,\\
							\{H^{*}(n)\}, & {\rm if}\quad 6\nmid n \quad {\rm but}\quad  4\vert n, \\
							\{(\frac{1}{3} \lceil \frac{n}{2} \rceil K_{3})\vee I_{\lfloor \frac{n}{2} \rfloor }\}, & 4\nmid n \quad {\rm but} \quad 3\vert \lceil \frac{n}{2} \rceil, \\
							\{(k_{1}K_{3}\cup k_{1}^{1}P_{2}\cup k_{1}^{2}P_{1})\vee I_{\lfloor \frac{n}{2} \rfloor}, \,H^{*}(n),\, (S_{k_{2}^{'}}\cup k_{2}K_{3})\vee I_{\lceil \frac{n-1}{2} \rceil}\}, & {\rm otherwise},	
						\end{array}
						\right.
					\end{split}
				\end{align*}
				where $1\leq \lceil \frac{n}{2}\rceil -3k_1\leq 2$, $k_1^1=\lceil \frac{n}{2}\rceil -3k_1-1$, $k_1^2=3k_1+2-\lceil \frac{n}{2}\rceil$, $0\leq 3k_2\leq \lceil \frac{n}{2}\rceil$ and $3k_2+1+k'_2=\lceil \frac{n}{2}\rceil$.
			\end{enumerate}
		\end{theorem}

	\section{Preliminaries}
	The key idea of our proofs is using a result of Simonovits (Theorem 2.1) to get a good vertex partition of an extremal graph for $C_{k,\;\ell}^{p+1}$. This was recommended by Liu in \cite{ref8}. The vertex split graphs family and the decomposition family result was introduced by Liu, which is also very crucial in the proofs. 
	
	Given a graph $H$ and a vertex $ v\in V(H) $ with $ d_H(v)\geq 2$, a vertex split on the vertex $v$ is defined as follows:
	replace  $v$ with an independent set of size $\lvert N_H(v)\rvert $ and each vertex is adjacent to exactly one distinct vertex in $N_H(v)$. Let $U$ be a vertex subset $U\subseteq  V(H)$, a vertex split on $U$ means applying vertex split on the vertices in $U$ one by one. Appearantly, the order of vertices we apply vertex split does not matter. $\mathcal{H}(H)$ is defined as the family of all the graphs which can be obta{\large }ined by applying vertex split on any vertex subset  $U\subseteq V(H)$. It is easy to see that $U$ can be empty, therefore, $H\in \mathcal{H}(H)$. $\mathcal {H}_p(H) $ is defined as the family of all the graphs obtained from $H$ by applying vertex split on the vertex subset $U \subseteq V(H)$, which satifies $ \chi(H[U])\leq p $. In particularly, $\mathcal{H}^*(H)$ is the family of all the graphs obtained from $H$ by applying vertex split on any independent set of $H$. It is not difficult to see that when $p\geq 2$ we have $\mathcal{H}^{*}(H) \subseteq \mathcal{H}_{p-1}(H)\subseteq \mathcal{H}(H) $.
	\begin{definition} (\cite{ref13})
			Given a family $\mathcal{L}$, define $ p=p(\mathcal{L})=\min\limits_{L\in \mathcal{L}}\chi(L)-1 $. Let $\mathcal{M}:=\mathcal{M}(\mathcal{L})$ be the family of minimal graphs $M$ for which there exist an $L \in \mathcal{L}$ and $t=t(L)$ satisfy that $L \subseteq {M^{'}} \vee {K_{p-1}}(t,t, \ldots ,t) $ where ${M^{'}}=M \cup I_{t}$. We call this the decomposition family of $\mathcal{L}$.
	\end{definition}

The following characterization of $\mathcal{M}(H^{p+1})$ was provided in \cite{ref10}.

\begin{lemma}(\cite{ref10})
	Let $H$ be any graph and $p \geq 2$ be any integer.
	\begin{enumerate}[(i)]
		\item If $p \geq 3$ and $\chi (H) \leq p-1$, then $\mathcal{M}(H^{p+1})=\mathcal{H}(H)$.
		\item If $p \geq 3$ and $\chi (H) = p$, then $\mathcal{M}(H^{p+1})=\mathcal{H}_{p-1}(H)$.
		\item If $p = 2$ and $\chi (H) = 2$ or $3$, then $\mathcal{M}(H^{p+1})=\mathcal{H}^{*}(H)$.
	\end{enumerate}
\end{lemma} 

	\begin{definition} 
	\cite{ref12} Let $G$ be a graph, $H_1$, $H_2$ be subgraphs of $G$. They are called symmetric if $H_1=H_2$ or satisfy the following conditions:
		\begin{enumerate}[(i)]
			\item $V(H_1)\cap V(H_2)=\emptyset $.
			\item $xy\notin E(G)$ if $x\in V(H_1)$, $y\in V(H_2)$.
			\item There exists an isomorphism $\omega $: $H_1\to H_2 $ such that for every $x\in V(H_1)$ and $u\in G-V(H_1)-V(H_2)$, $x$ is adjacent to $u$ if and only if $\omega (x)$ is adjacent to $u$.
		\end{enumerate}
	\end{definition}
	
	\begin{definition}(\cite{ref12})
		Let $\mathbb{D}(n,\,p,\,r)$ be the family of  graph $G$ of order $n$ with following conditions:
		\begin{enumerate}[(i)]
			\item It is possible to omit at most $r$ vertices of $G$ so that the remaining graph $G^{'}$ is a product of almost equal order: $G^{'}=\bigvee _{i\leq p} G^{i}$, where $\lvert V(G^{i})\rvert =n_{i}$ and $\lvert n_{i} -\frac{n}{p}\rvert \leq r$ $(1\leq i\leq p)$.
			\item For every $1\leq i\leq p$, there exist connected graphs $H_{i}$ such that $G^{i}=k_{i}H_{i}$, where $k_{i}=\frac{n_i}{\lvert V(H_i)\rvert }$ and any two copies $H_i^j$, $H_i^l$ in $G^{i}$ $(1\leq j<l\leq k_i)$ are symmetric subgraphs of $G$.
		
		\end{enumerate}
	\end{definition}
	
		The graphs $H_i$ in Definition 2.3 will be called the  blocks, the vertices in $G-G^{'}$ will be called exceptional vertices. Let $A_1,\cdots,A_p $ be the $p$ classes in $G'$, where $A_i=V(G^i)$ for any $i\in \{1,\cdots, p\}$. Let $W$ be the set of vertices in $G-G'$ that are adjacent to all vertices in $G'$ and let $B_{i}$  be the set of vertices in $G-G'-W$ that are not adjacent to any  vertex in $A_{i}$. 

	
	
		
	\begin{theorem}(\cite{ref12})
		Assume that a finite family $\mathcal{L}$ of forbidden graphs with $p(\mathcal{L})=p$ is given. If there exists some $ L\in \mathcal{L}$ with $m:=\lvert V(L) \rvert $ such that
		\begin{equation}
			L\subseteq P_{m}\vee K_{p-1} (m,\,m,\,\cdots ,m),
		\end{equation}
		then there exist $r=r(L)$ and $n_{0}=n_{0}(r)$ such that $\mathbb{D}(n,\,p,\,r)$ contains an $\mathcal{L}$-extremal graph for every $n>n_{0}$. Furthermore, if this is the only extremal graph in $\mathbb{D}(n,\,p,\,r)$, then it is the unique extremal graph for every sufficiently large n.
	\end{theorem}

		Let $\mathcal{Y}_{k+1,\;\ell+1}$ be the family of graphs obtained from a path $P_{k+1}$ by appending a path $P_{\ell +1}$ to one of its vertices except the end points. The vertex of degree 3 is called branching vertex.
	\begin{lemma}
		 For any integers $ k\geq 3$, $\ell \geq 0$, $ p \geq 2$ and $m= \lvert V(C_{k,\;\ell}^{p+1})\rvert$, 
		 
		 {\rm (i)}  if $Y$ is a graph in $\mathcal{Y}_{k+1,\;\ell+1}$, then we have $C_{k,\;\ell}^{p+1} \subseteq Y \vee K_{p-1}(m,m,\cdots ,m)$;
		 
		 {\rm (ii)}  we have $C_{k,\;\ell}^{p+1} \subseteq (P_{k+1}\cup P_{\ell +1}) \vee K_{p-1}(m,m,\cdots ,m)$.
	\end{lemma}
	\begin{proof1}
The fact $\chi (C_{k,\;\ell})=2$ or $3$, and Lemma 2.1 (ii) and (iii) imply that $\mathcal{H}^{*}(C_{k,\;\ell})\subseteq \mathcal{H}_{p-1}(C_{k,\;\ell})=\mathcal{M}(C_{k,\;\ell}^{p+1})$.

(i) 	Note that $Y$ can be obtained by applying vertex split on the one of the vertices except the center vertex of the cycle of $C_{k,\;\ell}$. Thus, $Y \in \mathcal{H}^{*}(C_{k,\;\ell})$  and then $Y\in \mathcal{M}(C_{k,\;\ell}^{p+1})$. By the definition of $\mathcal{M}(C_{k,\;\ell}^{p+1})$, we have $C_{k,\;\ell}^{p+1}\subseteq Y\vee K_{p-1}(m,m,\cdots,m)$.

(ii)  By applying vertex split on the center of $C_{k,\;\ell}$, the resulting  graph is $P_{k+1}\cup P_{\ell +1}$. Thus, $(P_{k+1}\cup P_{\ell +1})\in \mathcal{H}^{*}(C_{k,\;\ell})$ and then $(P_{k+1}\cup P_{\ell +1})\in \mathcal{M}(C_{k,\;\ell}^{p+1})$. By the definition of $\mathcal{M}(C_{k,\;\ell}^{p+1})$, we have $C_{k,\;\ell}^{p+1}\subseteq (P_{k+1}\cup P_{\ell +1})\vee K_{p-1}(m,m,\cdots,m)$.
	\end{proof1}

 Obviously $(P_{k+1}\cup P_{\ell +1})\subseteq P_m$, hence $C_{k,\;\ell}^{p+1} \subseteq P_{m} \vee K_{p-1}(m,m,\cdots ,m)$. By Theorem 2.1 we have ${\rm EX}(n,\,C_{k,\;\ell}^{p+1})\in \mathbb{D}(n,\,p,\,r)$.


	\section{Main Results}
	Let $G$ be an extremal graph for $C_{k,\;\ell}^{p+1}$, then we have $G\in \mathbb{D}(n,\,p,\,r)$. In the rest part, we always let $A_i$, $H_i$, $B_i$, $W$ be the decompositions of $G$ as defined in Definition 2.3. In this section we mainly further characterize the structure of them to find the extremal graph for $C_{k,\;\ell}^{p+1}$.
    Since $G\in \mathbb{D}(n,\,p,\,r)$, we may have the following upper bound for $e(G)$.
	\begin{lemma}
	Suppose $p\geq 2$, $n$ is sufficiently large, each block $H_i$ in $G^i$ is a single vertex, then we have
	\begin{equation}
		e(G)\leq e(T_p(n))+\frac{n\lvert W\rvert }{p}+o(n).
	\end{equation}
	\end{lemma}

    The following Lemma 3.2 implies a lower bound for $e(G)$ (see Corollary 3.1). We always write $t=\lfloor \frac{k-1}{2}\rfloor +\lfloor \frac{\ell-1}{2}\rfloor $, $m=\lvert V(C_{k,\;\ell}^{p+1})\rvert $ in this section.
	\begin{lemma}
		Suppose $k\geq 3$, $\ell\geq 2$, $p\geq 2$, n is sufficiently large. 
		\begin{enumerate}[(i)]
			\item $H(n,\,p,\,t+1)$ is $C_{k,\;\ell}^{p+1}$-free when $k$ is odd,  $\ell$ is odd.
			\item $H(n,\,p,\,t+2)$ is $C_{k,\;\ell}^{p+1}$-free when $k$ is odd,  $\ell$ is even.
			\item $H^{'}(n,\,p,\,t+1)$ is $C_{k,\;\ell}^{p+1}$-free when $k$ is even,  $\ell$ is odd.
			\item $H(n,\,p,\,t+2)$ is $C_{k,\;\ell}^{p+1}$-free when $k$ is even,  $\ell$ is even.
		\end{enumerate}
	\end{lemma}
	\begin{proof1}
		Denote by $Q$ the vertex set $V(K_{q-1})$ in graph $H(n,\,p,\,q)$.
		Any $(p+1)$-clique in $H(n,\,p,\,q)$ contains at least one vertex of $Q$ and there is at most one $(p+1)$-clique in $H'(n,\,p,\,q)$ has no vertex in $Q$. On the other hand, in $C_{k,\;\ell}^{p+1}$ there are only three ($p+1$)-cliques  which share one vertex, and any other pairs of ($p+1$)-cliques share at most one vertex.
		\begin{enumerate}[(i)]
			\item When $k$ is odd, $\ell$ is odd, we have $t=\frac{k+\ell-2}{2}  $. Note that $e(C_{k,\;\ell})=k+\ell$. In $H(n,\,p,\,t+1)$, we have $\lvert Q\rvert =t$, and then the number of $(p+1)$-cliques of $C_{k,\;\ell}^{p+1}$ in $H(n,\,p,\,t+1)$ is at most $2(t-1)+3=k+\ell-1<k+\ell$.  Hence, $H(n,\,p,\,t+1) $ is $C_{k,\;\ell}^{p+1}$-free.
			\item When $k$ is odd, $\ell$ is even, we have $t=\frac{k+\ell-3}{2}$. If $C_{k,\;\ell}^{p+1}\subseteq H(n,\,p,\,t+2)$, then the vertex set $Q$ is a vertex cover of $C_{k,\;\ell}$. The minimum size of vertex cover of $C_{k,\;\ell}$ is $\frac{k+\ell+1}{2}$, while we have $\lvert Q\rvert =t+1=\frac{k+\ell-1}{2} <\frac{k+\ell+1}{2}$. Hence, $H(n,\,p,\,t+2) $ is $C_{k,\;\ell}^{p+1}$-free.
			\item When $k$ is even, $\ell$ is odd, we have $t=\frac{k+\ell-3}{2} $. In $H^{'}(n,\,p,\,t+1)$, the number of $(p+1)$-cliques of $C_{k,\;\ell}^{p+1}$ is at most $2(t-1)+3+1=k+\ell-1<k+\ell $. Hence, $H^{'}(n,\,p,\,t+1) $ is $C_{k,\;\ell}^{p+1}$-free.
			\item When $k$ is even, $\ell$ is even, we have $t=\frac{k+\ell}{2} -2$. In $H(n,\,p,\,t+2)$ the number of $(p+1)$-cliques of $C_{k,\;\ell}^{p+1}$ is at most $2t+3=k+\ell-1<k+\ell$. Hence, $H(n,\,p,\,t+2)$ is $C_{k,\;\ell}^{p+1}$-free.
		\end{enumerate}
\end{proof1}

	\begin{corollary}
		Suppose $p\geq 2$, $n$ is sufficiently large.
		\begin{enumerate}[(i)]
	   \item  If $\ell$ is odd, then we have 
		\begin{equation}
			e(G)\geq e(T_p(n))+\frac{tn}{p}+o(n).
		\end{equation}
	\item If $\ell$ is even, then we have 
		\begin{equation}
			e(G)\geq e(T_p(n))+\frac{(t+1)n}{p}+o(n).
		\end{equation}
	\end{enumerate}
	\end{corollary}

If each block $H_i$ in $G^i$ ($1\leq i\leq p$) is a single vertex, then we may give some characterizations for the set $W$ and $B_{i}$ (see Lemma 3.3 and Lemma 3.4) of the extremal graph $G$.

	\begin{lemma}
	Suppose $k\geq 3$, $\ell \geq 2$, $p\geq 2$, $n$ is sufficiently large and each block $H_i$ in $G^i$ ($1\leq i\leq p$) is a single vertex. Then we have 
	
	{\rm(i)} $\lvert W\rvert =t$ when $\ell$ is odd; 
	
	{\rm(ii)} $\lvert W\rvert =t+1$ when $\ell$ is even.
\end{lemma}

\begin{proof1}
	(i) When $\ell$ is odd, by (3.1) and (3.2), it is easy to see $\lvert W\rvert \geq t$. Furthermore, if $\lvert W\rvert \geq t+1$, then we may suppose $\{w_1,w_2,\cdots ,w_{t+1}\}\subseteq W$. The assumption that $n$ is sufficiently large ensures $\lvert A_i \rvert $ is sufficiently large, and then we may suppose $\{u_1,u_2,\cdots ,u_{t+3}\}\subseteq A_1$.
	When $k$ is odd, we may have a graph $Y\in \mathcal{Y}_{k+1,\;\ell+1}$ with $w_{\frac{k-1}{2}}$ as the branching vertex and 
	\begin{equation*}
		P_{k+1}=u_{1}w_{1}u_{2}w_{2} \cdots w_{\frac{k-1}{2}}u_{ \frac{k+1}{2} }w_{ \frac{k+1}{2} },
	\end{equation*}
\begin{equation*}
	P_{\ell+1}=w_{\frac{k-1}{2}}u_{ \frac{k+3}{2} }w_{\frac{k+3}{2}}u_{ \frac{k+5}{2}} \cdots w_{t+1}u_{t+2}.
\end{equation*}
When $k$ is even, we may have a graph $Y\in \mathcal{Y}_{k+1,\;\ell+1}$ with $w_{ \frac{k}{2}}$ as the branching vertex and 
\begin{equation*}
	P_{k+1}=u_{1}w_{1}u_{2}w_{2} \cdots u_{ \frac{k}{2}}w_{ \frac{k}{2}}u_{t+3 },
\end{equation*}
\begin{equation*}
	P_{\ell+1}=w_{\frac{k}{2}}u_{ \frac{k+2}{2} }w_{\frac{k+2}{2}}u_{ \frac{k+4}{2}} \cdots w_{t+1}u_{t+2}.
\end{equation*}
Then $Y\subseteq G[A_1\cup W]$. Furthermore $Y\vee K_{p-1}(m,\,m,\cdots ,m)\subseteq G$. On the other hand, by Lemma 2.2 (i) $C_{k,\;\ell}^{p+1}\subseteq Y\vee K_{p-1}(m,\,m,\cdots ,m)$. Therefore, we have a $C_{k,\;\ell}^{p+1}$ in $G$ and this contradiction shows $\lvert W\rvert = t$. 
	
	(ii) When $\ell$ is even, by (3.1) and (3.3), it is easy to see $\lvert W\rvert \geq t+1$. Furthermore, if $\lvert W\rvert \geq t+2$, we may suppose $\{w_1,w_2, \cdots, w_{t+2}\}\subseteq W$, $\{u_1,u_2 ,\cdots, u_{t+4}\}\subseteq A_1$. 
	When $k$ is odd, we have two paths in $G$ with
	\begin{equation*}
	P_{k+1}=u_{1}w_{1}u_{2}w_{2} \cdots u_{ \frac{k+1}{2}}w_{ \frac{k+1}{2}},
	\end{equation*}
	\begin{equation*}
		P_{\ell+1}=u_{\frac{k+3}{2}}w_{\frac{k+3}{2}} \cdots w_{t+2}u_{t+3}.
	\end{equation*}
When $k$ is even, we have two paths in $G$ with
	\begin{equation*}
	P_{k+1}=u_{1}w_{1}u_{2}w_{2} \cdots u_{ \frac{k}{2} }w_{ \frac{k}{2}}u_{t+4},
\end{equation*}
\begin{equation*}
	P_{\ell+1}=u_{\frac{k+2}{2}}w_{\frac{k+2}{2}} \cdots w_{t+2}u_{t+3}.
\end{equation*}
	Then $(P_{k+1}\cup P_{\ell +1})\subseteq G[A_1\cup W]$. Furthermore $(P_{k+1}\cup P_{\ell +1})\vee K_{p-1}(m,\,m,\,\cdots , m)\subseteq G$. By Lemma 2.2 (ii) we may obtain a $C_{k,\;\ell}^{p+1}$ in $G$ and this contradiction shows $\lvert W\rvert = t+1$. 
\end{proof1}

	\begin{lemma}
	Suppose $k\geq 3$, $\ell \geq 2$, $p\geq 2$, $n$ is sufficiently large and  each block $H_i$ in $G^i$ ($1\leq i\leq p$) is a single vertex. Then each vertex in $B_{i}$ is adjacent to all the vertices in $V(G^{'})\backslash A_{i}$.
	\end{lemma}
\begin{proof1}
	By Lemma 3.3, we have $\lvert W\rvert =t$ when $\ell$ is odd, $\lvert W\rvert =t+1$ when $\ell$ is even. Without loss of generality, we may suppose to the contrary that there is a vertex $v$ in $B_2$ which is not adjacent to some vertex in $A_1$. Since the blocks in $G^i$ are symmetric, then $v$ is not adjacent to any vertex in $A_1$. Then when $\ell$ is odd, we have
	\begin{equation*}
		\begin{split}
			e(G)\leq e(T_p(n))+\frac{nt}{p}-\lvert A_1\rvert +o(n)\\
			=e(T_p(n))+\frac{n(t-1)}{p} +o(n),
		\end{split}
	\end{equation*}
	which contradicts to (3.2).
	
		When $\ell$ is even
	\begin{equation*}
		\begin{split}
			e(G)\leq e(T_p(n))+\frac{n(t+1)}{p}-\lvert A_1\rvert  +o(n)\\
			=e(T_p(n))+\frac{nt}{p} +o(n),
		\end{split}
	\end{equation*}
	which contradicts to (3.3).
\end{proof1}
\begin{lemma}
	Suppose $k\geq 3$ is odd, $\ell\geq 2$, $p\geq 2$, $n$ is sufficiently large. If the extremal graph $G$ for $C_{k,\;\ell}^{p+1}$ satisfies the following conditions:
	\begin{enumerate}[(i)]
		\item when $\ell$ is odd, $\lvert W\rvert=t$; when $\ell$ is even, $\lvert W\rvert=t+1$,
		\item each vertex in $B_i$ is adjacent to all the vertices in $V(G')\backslash A_i$, 
		
		then we have $e(G[B_i])=0$ ($1\leq i\leq p$).
		\end{enumerate}
	\end{lemma}
	\begin{proof1}
		Suppose $W=\{w_1,\cdots ,\,w_t\}$ ($W=\{w_1,\cdots ,\,w_{t+1}\}$  resp.) when $\ell$ is odd ($\ell$ is even resp.) and $\{u_1,\cdots,u_{t+1}\}\subseteq A_1$. We first show that $G[B_i]$ is $P_4$-free. Suppose not, without loss of generality, let $P_4=x_1y_1x_2y_2 \subseteq G[B_2]$, then when $\ell$ is odd, we have a lollipop with $w_{\frac{k-3}{2}}$ as center vertex and 
		\begin{equation*}
			C_k=y_1u_1w_1u_2w_2\cdots w_{\frac{k-3}{2}}u_{\frac{k-1}{2}}x_2y_1,
		\end{equation*}
	\begin{equation*}
		P_{\ell +1}=w_{\frac{k-3}{2}}u_{\frac{k+1}{2}}\cdots w_{t-1}u_{t+1}.
	\end{equation*}
		When $\ell$ is even, we have a lollipop with $w_{\frac{k-3}{2}}$ as center vertex and 
		\begin{equation*}
			C_k=y_1u_1w_1u_2w_2\cdots w_{\frac{k-3}{2}}u_{\frac{k-1}{2}}x_2y_1,
		\end{equation*}
	\begin{equation*}
		P_{\ell +1}=w_{\frac{k-3}{2}}u_{\frac{k+1}{2}}\cdots w_{t-1}u_{t+1}w_t.
	\end{equation*}
	Each edge between $W$ and $A_1$ can be blown up into a $(p+1)$-clique by using vertices in $A_i$ ($2\leq i\leq p$). The edge $y_1u_1$ can be blown up by using vertex $x_1$ and vertices in $A_i$ ($3\leq i\leq p$) and the edge $u_{\frac{k-1}{2}}x_2$ can be blown up by using vertex $y_2$ and vertices in $A_i$ ($3\leq i\leq p$). The edge $x_2y_1$ can be blown up by using a vertex in $A_1$ and using vertices 
	in $A_i$ ($3\leq i\leq p$). Then there is a $C_{k,\;\ell}^{p+1}$ in $G$. Thus $G[B_i]$ is $P_4$-free.
	
		 If $e(G[B_i])\neq 0$, then there is an edge $e=xy\in E(G[B_i])$ which satisfies that there exists a vertex $w\in W$ with $wx\in E(G)$ or $wy\in E(G)$. Otherwise, let $G_1$ be the graph obtained from $G$ by deleting all edges of $\cup _{i\leq p}G[B_i]$ and adding all missing edges between $W$ and $\cup _{i\leq p}B_i $, then $G_1$ is $C_{k,\;\ell}^{p+1}$-free. In fact if there is a $C_{k,\;\ell}^{p+1}$ in $G_1$ then we may choose some vertices in $\cup _{i\leq p} A_i$ as the substitutes for the vertices in $\cup _{i\leq p} B_i$ to obtain a $C_{k,\;\ell}^{p+1}$ in $G$. So $G_{1}$ is $C_{k,\;\ell}^{p+1}$-free. Since $G[B_i]$ is $P_4$-free and by Gallai Theorem (see \cite{ref2}), we have $e(G[B_i])\leq \lvert B_i\rvert$. Furthermore, $\lvert W \rvert\geq 2$, hence
		\begin{equation*}
			e(G_1)\geq e(G)-\sum_{i\leq p}e(G[B_i])+2\sum_{i\leq p}\lvert B_i\rvert>e(G),
		\end{equation*}
		which is a contradiction to the definition of $G$.

		So if there is an edge $xy$ in $G[B_1]$, then we may suppose $w_{ \frac{k-1}{2}}$ is adjacent to $x$ and  find two paths. In fact when $\ell$ is odd, we may have a graph $Y\in \mathcal{Y}_{k+1,\;\ell+1}$ with $w_{ \frac{k-1}{2} }$ as branching vertex and 
			\begin{equation*}
			P_{k+1}=u_{1}w_{1}u_{2}w_{2} \cdots u_{ \frac{k-1}{2}}w_{ \frac{k-1}{2} }xy,
		\end{equation*}
		\begin{equation*}
			P_{\ell +1}=w_{ \frac{k-1}{2}}u_{\frac{k+1}{2} }\cdots w_{t}u_{t+1}.
			\end{equation*}
		When $\ell$ is even,we may have a graph $Y\in \mathcal{Y}_{k+1,\;\ell+1}$ with $w_{ \frac{k-1}{2} }$ as branching vertex and 
			\begin{equation*}
			P_{k+1}=u_{1}w_{1}u_{2}w_{2} \cdots u_{ \frac{k-1}{2}}w_{ \frac{k-1}{2} }xy,
		\end{equation*}
		\begin{equation*}
			P_{\ell+1}=w_{ \frac{k-1}{2}}u_{\frac{k+1}{2}}\cdots w_{t}u_{t+1}w_{t+1}.
			\end{equation*}
	Then $Y\subseteq G[W\cup A_{1}\cup B_{1}]$. Furthermore, $Y\vee K_{p-1}(m,\,m,\cdots,m)\subseteq G$. Therefore by Lemma 2.2 (i) we have $C_{k,\;\ell}^{p+1}\subseteq G$ which is a contradiction. So we have $e(G[B_{i}])=0$ ($1\leq i\leq p$).
	\end{proof1}

Let the set of vertices of the cycle of $C_{k,\;\ell}$ be $\{a_1,\,a_2,\cdots,a_{k}\}$, the set of vertices of the path of $C_{k,\;\ell}$ be $\{a_1,\,b_2,\cdots,b_{\ell +1}\}$.
	\begin{theorem}
		When $k\geq 3$ is odd, $\ell\geq 2$, $n$ is sufficiently large, $H(n,\,3,\,t+1)$ $(H(n,\,3,\,t+2)\; resp.)$ is the unique extremal graph for $C_{k,\;\ell}^{4}$ when $\ell$ is odd (even resp.).
	\end{theorem}
	\begin{proof1}	
		Suppose $\ell$ is odd, then $t=\frac{k+\ell-2}{2}$. Now we will prove each block $H_i$ is a single vertex ($i=1,\,2,\,3$). 
		
		If $P_{3}\subseteq H_{1}$, then $(\frac{k-1}{2}P_{3}\cup (\ell +1)P_{2})\subseteq G[A_{1}]$, and $(\frac{k-1}{2}P_{3}\cup (\ell +1)P_{2})\vee K_2(m,\,m)\subseteq G$. On the other hand, by applying vertex split on the vertex set $U=\{a_1,\,a_3,\cdots,a_k,\,b_2,\,b_3,\cdots,b_{\ell}\}$, the resulting graph is $\frac{k-1}{2}P_{3}\cup (\ell +1)P_{2}$, since $\chi(C_{k,\;\ell}[U])=2$, we have $(\frac{k-1}{2}P_{3}\cup (\ell +1)P_{2}) \in \mathcal{H}_{2}(C_{k,\;\ell})$. By the fact $\chi (C_{k,\;\ell})=3$ and Lemma 2.1 (ii), we have $\mathcal{H}_{2}(C_{k,\;\ell})=\mathcal {M}(C_{k,\;\ell}^{4})$, so $(\frac{k-1}{2}P_{3}\cup (\ell +1)P_{2})\in \mathcal {M}(C_{k,\;\ell}^{4})$. By using the definition of $\mathcal{M}$, we have $C_{k,\;\ell}^{4}\subseteq (\frac{k-1}{2}P_{3}\cup (\ell +1)P_{2})\vee K_{2}(m,\,m)\subseteq G$. This contradiction implies that $H_i$ is $P_3$-free.
		
		
		Now suppose $H_{1}= H_{2} = P_{2}$ and then $(k+\ell )P_2\vee ((k+\ell )P_2\vee I_{k+\ell})\subseteq G$. Note that for any graph $F$, $F^{4}\subseteq e(F)P_2\vee F$ holds. The fact $C_{k,\;\ell}\subseteq ((k+\ell)P_{2}\vee I_{k+\ell})$ implies that $C_{k,\;\ell}^{4}\subseteq (k+\ell)P_{2}\vee ((k+\ell)P_{2}\vee I_{k+\ell})\subseteq G$. 
		
		If $H_{1}= P_{2}$, $H_{2}= P_{1}$, $H_{3}= P_{1}$, then 
		\begin{equation*}
			\begin{split}
				e(G)\leq e(T_3(n))+\frac{n\lvert W \rvert}{3}+\frac{\lvert A_{1}\rvert}{2}+o(n)\\
				\leq e(T_3(n))+\frac{n\lvert W \rvert}{3}+\frac{n}{6}+o(n).
			\end{split}
		\end{equation*}
		From Corollary 3.1, we have $e(G)\geq e(T_3(n))+\frac{tn}{3}+o(n)$, so $\lvert W\rvert \geq t$. Indeed, if $\lvert W\rvert \geq t$, then 
		$(\,\frac{k-1}{2} P_{3}\,\cup \,(\ell+1)P_{2}\,)\subseteq G[W\,\cup\, A_{1}]$, and then $(\frac{k-1}{2}P_3\cup (\ell +1)P_2)\vee K_2(m,\,m)\subseteq G$. So we have $C_{k,\;\ell}^{4}\subseteq (\,\frac{k-1}{2} P_{3}\,\cup\, (\ell+1)P_{2})\vee K_{2}\,(m,\,m\,)\subseteq G$, which is a contradiction.
	
	Therefore we have each block $H_i$ is a single vertex ($i=1,2,3$), and then by Lemma 3.3, 3.4, 3.5, we have (i) $\lvert W\rvert=t$; (ii) each vertex in $B_i$ is adjacent to all the vertices in $V(G')\backslash A_i$; (iii) $e(G[B_i])=0$ for $i=1,\,2,\,3$.
	By the maximality of $G$ we have
		 $G= H(n,\,p,\,t+1)$. By using the similar arguments, when $\ell$ is even, we have $G=H(n,\,p,\,t+2)$.
	\end{proof1}

	\begin{lemma}
	When $k\geq 4$, $\ell \geq 2$, $p=2$,  each block $H_{i}$ in $G^i$ is a single vertex $( i=1,\, 2)$.
	\end{lemma}
	
	\begin{proof1}
		Suppose $\ell$ is odd, then when $k$ is odd, $t=\frac{k+\ell -2}{2}$; when $k$ is even, $t=\frac{k+\ell -3}{2}$.
\begin{enumerate}[(i)]
	\item We claim that $H_i$ is $P_3$-free. Suppose to the contrary that $P_3\subseteq H_1$. 
	
	When $k$ is even, we have $((t+1)P_3 \cup P_2)\subseteq G[A_1]$, and then $((t+1)P_3 \cup P_2)\vee I_m \subseteq G$. On the other hand, when $k$ is even we apply vertex split on the vertex set $U=\{a_1,\,a_3,\cdots,a_{k-1},\,b_3,\,b_5,\cdots,b_{\ell}\}$, the resulting graph is $(t+1)P_3 \cup P_2$, since $\chi(C_{k,\;\ell}[U])=1$, we have 
	$((t+1)P_3 \cup P_2)\in \mathcal{H}^{*}(C_{k,\;\ell})$. By the fact $\chi (C_{k,\;\ell})=2$ and Lemma 2.1 (iii), we have $\mathcal{M}(C_{k,\;\ell}^3)=\mathcal{H}^{*}(C_{k,\;\ell})$, so $((t+1)P_3 \cup P_2)\in \mathcal{M}(C_{k,\;\ell}^3)$. By using the definition of $\mathcal{M}$, we have $C_{k,\;\ell}^3\subseteq ((t+1)P_3 \cup P_2)\vee I_m\subseteq G$. So $H_1$ is $P_3$-free when $k$ is even.
	
	Now suppose $k$ is odd. If $P_3\subseteq H_1$ and $P_2\subseteq H_2$ then $(mP_3\vee mP_2)\subseteq G$ while $C_{k,\;\ell}^{3}\subseteq (mP_3\vee mP_2)$, so $H_2=P_1$.
	On the other hand, $(P_4\cup (\frac{k-3}{2} +\frac{\ell-1}{2})P_3\cup P_2)\in \mathcal{H}^{*}(C_{k,\;\ell})$. By Lemma 2.1 (iii) $\mathcal{H}^{*}(C_{k,\;\ell})=\mathcal{M}(C_{k,\;\ell}^3)$, we have $(P_4\cup (\frac{k-3}{2} +\frac{\ell-1}{2})P_3\cup P_2)\in \mathcal{M}(C_{k,\;\ell}^3)$. If $\lvert W\rvert \neq 0$, then we have $(P_4\cup (\frac{k-3}{2} +\frac{\ell-1}{2})P_3\cup P_2)\subseteq G[W\cup A_1]$ and $C_{k,\;\ell}^3\subseteq (P_4\cup (\frac{k-3}{2} +\frac{\ell-1}{2})P_3\cup P_2)\vee I_m\subseteq G$. If $P_4\subseteq H_1$, then we have $(P_4\cup (\frac{k-3}{2} +\frac{\ell-1}{2})P_3\cup P_2)\subseteq G[A_1]$ and $C_{k,\;\ell}^3\subseteq G$. Therefore, we have $\lvert W\rvert =0$ and $H_1$ is $P_4$-free.
	 Since $H_1$ is $P_4$-free, by Gallai Theorem (see [2]), the size of $G$ is maximized when $H_1=K_3$. Hence, 
	\begin{equation*}
		e(G)\leq e(T_2(n))+\lvert A_1\rvert +o(n)\leq e(T_2(n))+\frac{n}{2}+o(n).
	\end{equation*}
	While it contradicts (3.2) when $k\geq 4$ and $\ell \geq 2$. Thus $H_i$ is $P_3$-free.

	\item If $H_1=H_2=P_2$, then $\lvert W\rvert=0$ holds. Otherwise, let $w\in W$,   $u_iu'_i\subseteq G[A_1]$, $v_iv'_i\subseteq G[A_2]$ ($1\leq i\leq t+1$) and then we may find a $C_{k,\;\ell}^{3}$ in $G$. In fact when $k$ is odd, we have a lollipop with $v_{\frac{k-1}{2}}$ as center vertex and 
	\begin{equation*}
		C_k=wu_1v_1\cdots u_{\frac{k-1}{2}}v_{\frac{k-1}{2}}w,
		\end{equation*}
	\begin{equation*}
		P_{\ell +1}=v_{\frac{k-1}{2}}u_{\frac{k+1}{2}}\cdots v_tu_{t+1}.
	\end{equation*}
When $k$ is even, we have a lollipop with $u_{\frac{k}{2}}$ as center vertex and 
 \begin{equation*}
 	C_k=wu_1v_1\cdots u_{\frac{k}{2}}w,
 \end{equation*}
 \begin{equation*}
 	P_{\ell +1}=u_{\frac{k}{2}}v_{\frac{k}{2}}\cdots u_{t+1}v_{t+1}.
 \end{equation*}
The edges between $W$ and $A_1$ can be blown up into a triangle by using one vertex in $A_2$ and the edges  between $W$ and $A_2$ can be blown up by using one vertex in $A_1$. The edge $u_iv_i$ can be expanded by using $u'_i$ and the edge $v_iu_{i+1}$ can be expanded by using $v'_i$ ($1\leq i\leq t$). Therefore $\lvert W\rvert=0$. Furthermore, we have 
\begin{equation*}
	e(G)\leq e(T_2(n))+\frac{\lvert A_1\rvert}{2}+\frac{\lvert A_2\rvert}{2}+o(n)\\
    \leq e(T_2(n))+\frac{n}{2}+o(n),
\end{equation*}
while it contradicts (3.2).
\item If $H_1=P_2$, $H_2=P_1$, then 
\begin{equation*}
	\begin{split}
		e(G) \leq e(T_2(n))+\frac{n\lvert W \rvert}{2}+\frac{n}{4}+o(n).
	\end{split}
\end{equation*}
Recall Corollary 3.1, $e(G)\geq e(T_2(n))+\frac{tn}{2}+o(n)$, so we have 	$\lvert W \rvert \geq t$. On the other hand, if $\lvert W \rvert \geq t$, then $(P_{k+1}\cup P_{\ell+1})\subseteq G[W\cup A_{1}]$. The fact that $(P_{k+1}\cup P_{\ell+1})\vee I_m\subseteq G$ and Lemma 2.2 (ii) imply that $C_{k,\;\ell}^{p+1}\subseteq G$ which is a contradiction. Therefore we have $H_i=P_1$ ($i=1,\,2$) when $\ell$ is odd.
\end{enumerate}
By using the similar arguments, we may prove $H_i=P_1$ ($i=1,\,2$) when $\ell$ is even. 
\end{proof1}

	\begin{theorem}
		When $k\geq 4$, $\ell\geq 2$ is odd,  $H(n,\,2,\,t+1)$ $(H^{'}(n,\,2,\,t+1) \;resp.)$ is the unique extremal graph for $C_{k,\;\ell}^{3}$ when $k$ is odd (even resp.).
	\end{theorem}
\begin{proof1}
From Lemma 3.3, Lemma 3.4 and Lemma 3.6, we know that (i) $\lvert W\rvert=t$; (ii) each vertex in $B_i$ is adjacent to all vertices in $V(G')\backslash A_i$; (iii) each block $H_i$ is a single vertex ($i=1,2$). To find the extremal graph, we only need to characterize the subgraph $G[B_i]$. When $k$ is odd, by Lemma 3.5, we have $e(G[B_i])=0$ ($i=1,\,2$). When $k$ is even, we have the following claim.
	
		\noindent \textbf{Claim.}  $e(G[B_{1}])+e(G[B_{2}])= 1.$

		 	\noindent \textbf{Proof of Claim.}
	 Set $\{w_{1},\,\cdots, w_{t} \}=W$, $\{u_{1},\,\cdots,u_{t+2}\}\subseteq A_{1}$. We first show that $G[B_i]$ is $P_3$-free for $i=1,\,2$. Suppose not, without loss of generality, let $P_3=xyz\subseteq G[B_2]$, then we have a lollipop with $w_{\frac{k-2}{2}}$ as center vertex and
	 \begin{equation*}
	 	C_k=yu_1w_1\cdots w_{\frac{k-2}{2}}u_{\frac{k}{2}}y,
	 \end{equation*}
	 \begin{equation*}
	 	P_{\ell+1}=w_{\frac{k-2}{2}}u_{\frac{k+2}{2}}\cdots w_tu_{t+2}.
	 \end{equation*}
	 Ecah edge between $W$ and $A_1$ can be blown up into a triangle by using one vertex in $A_2$. The edge $yu_1$ can be blown up by using vertex $x$ and the edge $u_{\frac{k}{2}}y$ can be blown up by using vertex $z$. Then there is a $C_{k,\;\ell}^3$ in $G$. Thus $G[B_i]$ is $P_3$-free for $i=1,\,2$.
	 
	 Suppose to the contrary that there exist edges $e_{1},\;e_{2}$ in $G[B_{1}\cup B_{2}]$, say $e_{1}=x_{1}y_{1}$, $e_{2}=x_{2}y_{2}$. 
			\begin{enumerate}[(i)]
			\item Suppose that $\{e_{1}$, $e_{2}\}\subseteq E(G[B_1])$. Since $G[B_i]$ is $P_3$-free, then $e_1$ and $e_2$ are independent. Note that for any vertices $w_i,\,w_j\,\in W$, we have $\lvert \{x_1,\,y_1\}\cap\, (N(w_{i})\cup N(w_{j}))\rvert \geq 3$
			 Otherwise, if $\lvert \{x_1,\,y_1\}\cap\, (N(w_{i})\cup N(w_{j}))\rvert \leq 2$, we can obtain a graph $G_1$ by deleting  edge $e_1$ and adding missing edges between $\{w_i,\,w_j\}$ and $\{x_1,\,y_1\}$. $G_1$ is still $C_{k,\;\ell}^{3}$-free. In fact if there is a $C_{k,\;\ell}^{3}$ in $G_1$, since $G[B_1]$ is $P_3$-free, then we may choose two vertices in $A_1$ as the substitutes for $x_1,\,y_1$ to obtain a $C_{k,\;\ell}^3$ in $G$.  While $e(G_1)>e(G)$, which is a contradiction. Hence $\lvert \{x_1,\,y_1\}\cap\, (N(w_{i})\cup N(w_{j}))\rvert \geq 3$. Similarly, $\lvert \{x_2,\,y_2\}\cap\, (N(w_{i})\cup N(w_{j}))\rvert \geq 3$.
			We may suppose $w_1$ is adjacent to $y_1$ and $x_2$, $w_2$ is adjacent to $x_1$. Further, we have a graph $Y\in \mathcal{Y}_{k+1,\;\ell+1}$ with $w_{\frac{k-2}{2}}$ as branching vertex and 
			\begin{equation*}
				P_{k+1}=y_2x_2w_1y_1x_1w_2u_1\cdots w_{\frac{k-2}{2}}u_{\frac{k-4}{2}},
			\end{equation*}
			\begin{equation*}
				P_{\ell +1}=w_{\frac{k-2}{2}}u_{\frac{k-2}{2}}w_{\frac{k}{2}}\cdots w_tu_t,
			\end{equation*}
			then $Y\subseteq G[W\cup A_1\cup B_1]$ and then we have $Y\vee I_m \subseteq G$.  On the other hand, by Lemma 2.2 (i) we have $C_{k,\;\ell}^3\subseteq Y\vee I_m$, hence we have $C_{k,\;\ell}^{3}\subseteq G$.
			
			\item Suppose $e_{1}\in E(G[B_{1}])$ and $e_{2}\in E(G[B_{2}])$. Then there are at least three edges between $\{x_1,\,y_1\}$ and $\{x_2,\,y_2\}$. Otherwise the graph $G_{1}$ obtained by deleting edge $x_2y_2$ and adding the missing edges between $\{x_1,\,y_1\}$ and $\{x_2,\,y_2\}$ is still $C_{k,\;\ell}^3$-free. In fact if there is a $C_{k,\;\ell}^3$ in $G_1$ then we may choose two vertices in $A_2$ as the substitutes for $x_2$, $y_2$ to obtain a $C_{k,\;\ell}^3$ in $G$, while $e(G_1)>e(G)$. So we may suppose $x_2$ is adjacent to $x_1$ and $y_1$.
			Furthermore, as the proof  of Lemma 3.5, we may suppose that $w_{\frac{k-2}{2}}$ is adjacent to $x_1$ and we have a lollipop with $w_{\frac{k-2}{2}}$ as center vertex and 
			\begin{equation*}
				C_{k}=u_1w_1u_2w_2\cdots w_{\frac{k-2}{2}}x_1x_2u_1,
			\end{equation*}
			\begin{equation*}
				P_{\ell +1}=w_{\frac{k-2}{2}}u_{\frac{k}{2}}\cdots w_tu_{t+1}.
			\end{equation*}
			The edge between $W$ and $A_1\cup B_1$ can be expanded into a triangle by using one vertex in $A_2$. The edge $x_1x_2$ can be expanded by using vertex $y_1$ and the edge $x_2u_1$ can be expanded by using vertex $y_2$. Then there is a $C_{k,\;\ell}^3$ in $G$ which is a contradiction.
			\end{enumerate}
			 Thus $e(G[B_{1}])+e(G[B_{2}])= 1$. Therefore, when $k$ is odd, $G=H(n,\;2,\;t+1)$. When $k$ is even, $G=H'(n,\;2,\;t+1)$.
   \end{proof1}
	
	\begin{theorem}
		When $k\geq 4$, $\ell \geq 2$ is even, $H(n,\,2,\,t+2)$ is the unique extremal graph for $C_{k,\;\ell}^{3}$.
	\end{theorem}
	\begin{proof1}
	 From Lemma 3.3, Lemma 3.4 and Lemma 3.6, we know that (i) $\lvert W\rvert=t+1$; (ii) each vertex in $B_i$ is adjacent to all vertices in $V(G')\backslash A_i$; (iii) each block $H_i$ in $A_i$ is a single vertex ($i=1,\,2$). To finish the proof, we only need to prove $e(G[B_i])=0$, $i=1,\,2$. When $k$ is odd, by Lemma 3.5 we have $e(G[B_i])=0$. When $k$ is even, we have the following claim.
	 	\textbf{Claim.} $e(G[B_i])=0$,  $i=1,\,2$.
	 	
	 	\textbf{Proof of Claim.}  Set $\{w_1,\cdots,w_{t+1}\}=W$, $\{u_1,\cdots,u_{t+2}\}\subseteq A_1$. Firstly, $G[B_i]$ is $P_3$-free for $i=1,\,2$. Suppose not, without loss of generality, let $P_3=xyz\subseteq G[B_2]$, then we have a lollipop with $w_{\frac{k-2}{2}}$ center vertex and
	 	\begin{equation*}
	 		C_k=yu_1w_1\cdots w_{\frac{k-2}{2}}u_{\frac{k}{2}}y,
	 	\end{equation*}
	 		\begin{equation*}
	 			P_{\ell +1}=w_{\frac{k-2}{2}}u_{\frac{k+2}{2}}\cdots u_{t+2}w_{t+1}.
	 	\end{equation*}
	 The edge between $W$ and $A_1$ can be blown up into a triangle by using one vertex in $A_2$. The edge $yu_1$ can be blown up by using vertex $x$ and the edge $u_{\frac{k}{2}}y$ can be blown up by using vertex $z$. Then there is a $C_{k,\;\ell}^3$ in $G$. Thus $G[B_i]$ is $P_3$-free for $i=1,\,2$.  Suppose to the contrary that there is an edge $xy$ in $G[B_1]$. Then we may suppose there exists a vertex $w\in W$ with $wx\in e(G)$ or $wy\in e(G)$, otherwise we may have a graph $G_1$ obtained from $G$ by deleting all edges of $G[B_1]\cup G[B_2]$ and adding all missing edges between $W$ and $B_1\cup  B_2$, while $e(G_1)>e(G)$. Then we may suppose $w_{t+1}$ is adjacent to $x$ and have a graph $Y\in \mathcal{Y}_{k+1,\;\ell+1}$ with $w_{\frac{k}{2}}$ as branching vertex and 
				\begin{equation*}
				P_{k+1}=u_{1}w_{1}u_{2}w_{2}\cdots w_{\frac{k}{2}}u_{t+2},
			\end{equation*}
			\begin{equation*}
				P_{\ell+1}=w_{\frac{k}{2}}u_{\frac{k+2}{2}}w_{\frac{k+2}{2}}\cdots w_{t+1}xy.
			\end{equation*}
		Then $Y\subseteq G[W\cup A_1\cup B_1]$.
		Furthermore, $Y\vee I_m\subseteq G$. By Lemma 2.2 (i) we have $C_{k,\;\ell}^{p+1}\subseteq G$ which is a contradiction. So we have $e(G[B_i])=0$.
	
		By the maximality of $G$, we have $G= H(n,\,2,\,t+2)$.
	\end{proof1}
	
	\begin{theorem}
		When $\ell \geq 2$, $H(n,\,2,\,t+1)$ $(\,H(n,\,2,\,t+2) \; resp.)$ is the unique extremal graph for $C_{3,\;\ell}^{3}$ when $\ell$ is odd  (even resp.).
	\end{theorem}

	\begin{proof1}
		We first characterize the subgraph of $H_i$ and $G[B_i]$ ($i=1,\,2$).
		
		\textbf{Claim.} Each block $H_i$ is a single vertex ($i=1,\,2$).
	
	\textbf{Proof of Claim.}
		\begin{enumerate}[(i)]
		\item We first show that $H_{1}$ and $H_{2}$ are $P_{4}$-free. Without loss of generality, if $P_{4}\subseteq H_{1}$, then we have $P_{4}\,\cup \,\lfloor \frac{\ell}{3} \rfloor P_{4}\,\cup\, P_{\ell-3\lfloor \frac{\ell}{3} \rfloor +1}\subseteq G[A_{1}]$. By Lemma 2.1 (iii) and $\chi(C_{k,\;\ell})=3$, we have $\mathcal{M}(C_{k,\;\ell}^3)=\mathcal{H}^{*}(C_{k,\;\ell})$.
		Since $P_{4}\,\cup \,\lfloor \frac{\ell}{3} \rfloor P_{4}\,\cup\, P_{\ell-3\lfloor \frac{\ell}{3} \rfloor +1}\in \mathcal{H}^{*}(C_{3,\;\ell})$ and  $\mathcal{M}(C_{3,\;\ell }^{3})=\mathcal{H}^{*}(C_{3,\;\ell})$, we have  $C_{3,\;\ell}^{3}\subseteq G[A_1]\bigvee I_{m}\subseteq G$, which is a contradiction.
		\item  Secondly we show that $H_{1}$ and $H_{2}$ are $P_{3}$-free.
		Suppose to the contrary that $P_{3}\subseteq H_{1}$. We will prove  $H_{2}= P_1$ and $\lvert W \rvert =0$. Otherwise, if $P_{2}\subseteq H_{2}$, then we may let $u'_iu_iu''_i\subseteq G[A_1]$, $v_iv'_i\subseteq G[A_2]$ ($1\leq i\leq \ell$). Then we have $C_3=u_1v_1v'_1u_1$. When $\ell$ is odd, we have
				\begin{equation*}
					P_{\ell+1}=v_1u_2v_2\cdots u_{ \frac{\ell +3}{2}},
				\end{equation*}
			and	when $\ell$ is even, we have
				\begin{equation*}
					P_{\ell+1}=v_1u_2v_2\cdots u_{\frac{\ell+2}{2} }v_{ \frac{\ell+2}{2}}.
				\end{equation*}
				The edge $v_1v'_1$ may be expanded into a triangle by using one vertex in $A_1$; the edge $u_1v_1$ may be expanded by using vertex $u'_1$; the edge $u_1v'_1$ may be expanded by using vertex $u''_1$. The edge $v_iu_{i+1}$ may be expanded by using vertex $u'_i$ and the edge $u_iv_{i}$ may be expanded by using vertex $u''_i$. Then there is a $C_{3,\;\ell}^{3}\subseteq G[A_{1}\cup A_{2}]\subseteq G$. So we have $H_2=P_1$.
				
				If $\lvert W \rvert \neq 0$, then there is a $w\in W$. When $\ell$ is odd, it is easy to see $(P_4\cup \frac{\ell-1}{2}P_3\cup P_2)\subseteq G[W\cup A_1]$. On the other hand, $(P_4\cup \frac{\ell-1}{2}P_3\cup P_2)\in \mathcal{H}^{*}(C_{3,\;\ell})$. By Lemma 2.1 (iii) $\mathcal{H}^{*}(C_{3,\;\ell})=\mathcal{M}(C_{3,\;\ell}^3)$, so $(P_4\cup \frac{\ell-1}{2}P_3\cup P_2)\in \mathcal{M}(C_{3,\;\ell}^3)$. Then we have $C_{3,\;\ell}^3\subseteq (P_4\cup \frac{\ell-1}{2}P_3\cup P_2)\vee I_t\subseteq G$. By using the same arguments we may prove the results for even $\ell$.
			So we have $H_2=P_1$ and $\vert W\rvert =0$ and from (i) we know that $G^i$ is $P_4$-free.
				By Gallai Theorem (see \cite{ref2}), we have
				\begin{equation*}
				 e(G)\leq e(T_2(n))+\lvert A_1\rvert +o(n)\leq e(T_2(n))+\frac{n}{2}+o(n).
				 \end{equation*}
				  While it contradicts Corollary 3.1. So $H_1$ and $H_2$ are $P_3$-free.
				
				\item Now supose $H_{1}=P_{2}$. 
				We may set $u_iu'_i\subseteq G[A_1]$ ($1\leq i\leq \ell$). 
				 We distinguish the following two cases  according to $\ell$.
				
				\textbf{Case 1}: $\ell=2,\;3$. Set $\{v_1,v_2,\cdots,v_{\ell }\}\subseteq A_2$. We claim  $\lvert W \rvert =0$. If not, suppose $w\in W$. Let $C_3=wu_1u'_1w$, $P_4=wu_2v_2u_3$. Each edge of $C_3$ and the edge $wu_2$ can be blown up into a triangle by employing one vertex in $A_2$, the edge $u_2v_2$ can be blown up by using vertex $u'_2$, the edge $v_2u_3$ can be blown up by using vertex $u'_3$, then there is a $C_{3,\;\ell}^3$ in $G$. So $\lvert W\rvert =0$, and then  we have 
				\begin{equation*}
				e(G)\leq e(T_2(n))+\frac{\lvert A_1\rvert}{2} +\frac{\lvert A_2\rvert}{2}+o(n)\leq e(T_2(n))+\frac{n}{2}+o(n),
				\end{equation*}
				while it contradicts Corollary 3.1. So $H_1$ is a single vertex.
			
				
				\textbf{Case 2}: $\ell=4a+b$ ($a\geq 1$), we  may claim  $\lvert W\rvert\leq a$. If $\lvert W\rvert\geq a+1$, let  $\{w_{1},\,w_{2},\cdots, w_{a+1}\}\subseteq W$ and $\{v_{1},\,v_{2},\cdots ,v_{\ell }\}\subseteq A_{2}$, then we have a  $C_3=w_1u_1u'_1w_1$, each edge of $C_3$ can be blown up into a triangle.
				
				When $\ell =4a$, set
				\begin{equation*}
					P_{\ell +1}=w_1u_2v_2u_3w_2u_4v_4u_5\cdots w_{a+1}.
					\end{equation*}
				When $\ell =4a+1$, set
				\begin{equation*}
					P_{\ell +1}=w_1u_2v_2u_3w_2u_4v_4u_5\cdots w_{a+1}u_{2a+2}.
				\end{equation*}
				When $\ell =4a+2$, set
				\begin{equation*}
					P_{\ell +1}=w_1u_2v_2u_3w_2u_4v_4u_5\cdots w_{a+1}u_{2a+2}v_{2a+2}.
				\end{equation*}
				When $\ell =4a+3$, set
					\begin{equation*}
					P_{\ell +1}=w_1u_2v_2u_3w_2u_4v_4u_5\cdots w_{a+1}u_{2a+2}v_{2a+2}u_{2a+3}.
				\end{equation*}
			The edges $w_{i}u_{2i}$ and $u_{2i-1}w_{i}$ can be blown up into a triangle by using one vertex in $A_2$, the edge $u_iv_i$ can be blown up by using vertex $u'_i$, the edge $v_iu_{i+1}$ can be blown up by using vertex $u'_{i+1}$. Then there is a $C_{3,\;\ell}^3$ in $G$  and it comtradicts the definition of $G$. Therefore, we have $\lvert W\rvert \leq a$. 
			Next we will prove that if $\lvert W\rvert\neq 0$, then $H_2=P_1$. Otherwise, let $v_iv'_i \subseteq G[A_2]$ ($1\leq i\leq \ell$) and  $w\in W$, then we have a  $C_3=wu_1u'_1w$. When $\ell$ is odd, we have
			\begin{equation*}
				P_{\ell +1}=u_1v_1\cdots u_{\frac{\ell +1}{2}}v_{\frac{\ell +1}{2}},
			\end{equation*}
			when $\ell$ is even, we have
			\begin{equation*}
				P_{\ell +1}=u_1v_1\cdots u_{\frac{\ell}{2}}v_{\frac{\ell }{2}}u_{\frac{\ell +2}{2}}.
			\end{equation*}
			
			Each edge of $C_3$ can be blown up by employing one vertex in $A_2$, the edge $u_iv_i$ can be blown up by using vertex $v'_i$, and the edge $v_ju_{j+1}$ can be blown up by using vertex $u'_{j+1}$. Then there is a $C_{3,\;\ell}^3$ in $G$ which is a contradiction.
			
			When $\lvert W\rvert=0$, from (ii) we know that $H_i$ is $P_3$-free, so we have
			\begin{equation*}
				e(G)\leq e(T_2(n))+\frac{\lvert A_1\rvert}{2}+\frac{\lvert A_2\rvert}{2}+o(n)\\
				\leq e(T_2(n))+\frac{n(\lfloor \frac{\ell }{4}\rfloor+1)}{2}+o(n).
			\end{equation*}
		When $0<\lvert W\rvert \leq a$, we have
			\begin{equation*}
				e(G)\leq e(T_2(n))+\frac{n\lvert W\rvert}{2}+\frac{\lvert A_1\rvert}{2}+o(n)\\
				\leq e(T_2(n))+\frac{n(\lfloor \frac{\ell }{4}\rfloor+1)}{2}+o(n).
			\end{equation*}
		When $\ell \geq 4$ is odd, $\lfloor \frac{\ell }{4}\rfloor +1 <\frac{\ell +1}{2}$, when $\ell \geq 4$ is even, $\lfloor \frac{\ell }{4}\rfloor+1 <\frac{\ell}{2}+1$ and these contradict Corollary 3.1. Therefore $H_{1}=H_{2}=P_{1}$.
		\end{enumerate}
	   Now we have proved that $H_1=H_2=P_1$. Combining Lemma 3.3 - 3.5, we have $G= H(n,\,2,\,t+1)$ when $\ell$ is odd; $G= H(n,\,2,\,t+2)$ when $\ell$ is even.
		\end{proof1}
	

\begin{thebibliography}{99}
	\small{	\bibitem{ref1}G. Chen, R. J. Gould, F. Pfender, B. Wei, Extremal graphs for intersecting cliques, J. Comb. Theory, Ser. B 89 (2003) 159-171.
		\bibitem{ref2}P. Erd\H{o}s, T. Gallai, On maximal paths and circuits of graphs, Acta Math. Acad. Sci. Hung. 10, (1959) 337-356. 
		\bibitem{ref3}P. Erd\H{o}s, M. Simonovits, A limit theorem in graph theory, Studia Sci. Math. Hung. 1 (1966) 51-57.
		\bibitem{ref4}P. Erd\H{o}s, $\ddot{U}$ber ein Extremalproblem in der Graphentheorie (German), Arch. Math. (Basel) 13 (1962) 122-127.
		\bibitem{ref5}P. Erd\H{o}s, Z. F$\ddot{u}$redi, R. J. Gould, D. S. Gunderson, Extremal graphs for intersecting triangles, J. Comb. Theory, Ser. B 63 (1995) 89-100.
		\bibitem{ref6}R. Glebov, Extremal graphs for clique-paths, arXiv: 1111.7029v1, 2011.
		\bibitem{ref7}X. Hou, Y. Qiu, B. Liu, Extremal graph for intersecting odd cycles, Electron. J. Comb. 23 (2) (2016) P2. 29.
		\bibitem{ref8}H. Liu, Extremal graphs for blow-ups of cycles and trees, Electron. J. Comb. 20 (1) (2013) P65.
		\bibitem{ref9}J. W. Moon, On independent complete subgraphs in a graph, Can. J. Math. 20 (1968) 95-102.
		\bibitem{ref10}Z. Ni, L. Kang, E. Shan, Extremal graphs for blow-ups of keyrings, Graphs Comb. 36 (2020) 1827-1853.
		\bibitem{ref11}M. Simonovits, A method for solving extremal problems in graph theory, stablity problems, in: Theory of Graphs, Proc. Colloq., Tihany, 1996, Academic Press, New York, 1968, pp. 279-319.
		\bibitem{ref12}M. Simonovits, Extremal graph problems with symmetrical extremal graphs, additional chromatic conditions, Discrete Math. 7 (1974) 349-376.
		\bibitem{ref13}P. Tur\'{a}n, On the theory of graphs, Colloquium Mathematicum. 3 (1) (1954) 19-30.
		\bibitem{ref14}A. Wang, X. Hou, B. Liu, The Tur\'{a}n number for the edge blow-up of trees, Discrete Math. 344 (2021) 112627.
		\bibitem{ref15}L. Yuan, Extremal graphs for edge blow-up of graphs, J. Comb. Theory, Ser. B 152 (2022) 379-398.
		\bibitem{ref16}H. Zhu, L. Kang, E. Shan, Extremal graphs for odd-ballooning of paths and cycles, Graphs Comb. 36 (2020) 755-766.}
	\end{thebibliography}
\end{document}